\documentclass[12pt]{article}
\usepackage{stmaryrd}
\usepackage{bbm}
\usepackage[french,english]{babel}

\usepackage{setspace,mathtools,amsmath,amssymb,amsthm,fullpage,outline,centernot,slashed} 
\usepackage[pdftex]{graphicx}
\usepackage{mathrsfs}
\usepackage{enumerate}
\usepackage{accents}
\usepackage[colorlinks=true,citecolor=black,linkcolor=black,urlcolor=blue]{hyperref}
\usepackage{upgreek}

\theoremstyle{plain}
\newtheorem{theorem}{Theorem}
\newtheorem{lemma}[theorem]{Lemma}

\theoremstyle{definition}

\newtheorem{question}[theorem]{Question}
\theoremstyle{remark}
\newtheorem{remark}[theorem]{Remark}

\newcommand{\R}{\mathbb{R}}
\newcommand{\N}{\mathbb{N}}

\newcommand{\ben}{\begin{enumerate}}
\newcommand{\een}{\end{enumerate}}
\newcommand{\bit}{\begin{itemize}}
\newcommand{\eit}{\end{itemize}}
\renewcommand{\>}{\rangle} 
\def\bal#1\eal{\begin{align*}#1\end{align*}}

\DeclareMathSymbol{\mlq}{\mathord}{operators}{``}
\DeclareMathSymbol{\mrq}{\mathord}{operators}{`'}

\def\Xint#1{\mathchoice
{\XXint\displaystyle\textstyle{#1}}%
{\XXint\textstyle\scriptstyle{#1}}%
{\XXint\scriptstyle\scriptscriptstyle{#1}}%
{\XXint\scriptscriptstyle\scriptscriptstyle{#1}}%
\!\int}
\def\XXint#1#2#3{{\setbox0=\hbox{$#1{#2#3}{\int}$ }
\vcenter{\hbox{$#2#3$ }}\kern-.6\wd0}}

\def\dashint{\Xint-}

\title{Lower bounds and fixed points for the centered Hardy--Littlewood maximal operator}
\author{Samuel Zbarsky, Princeton University}

\begin{document}
\maketitle
\begin{abstract}
For all $p>1$ and all centrally symmetric convex bodies $K\subset \R^d$ define $Mf$ as the centered maximal function associated to $K$. We show that when $d=1$ or $d=2$, we have $||Mf||_p\ge (1+\epsilon(p,K))||f||_p$. For $d\ge 3$, let $q_0(K)$ be the infimum value of $p$ for which $M$ has a fixed point. We show that for generic shapes $K$, we have $q_0(K)>q_0(B(0,1))$.
\end{abstract}
\section{Introduction}
Let $K$ be some centrally symmetric convex body, and define the centered maximal function of a locally integrable function $f$ on $\R^d$ by
\begin{equation}\label{Mdef}
Mf(x)=\sup_{\lambda>0, S=x+\lambda K}\frac{1}{|S|}\int_S f.
\end{equation}
where $|S|$ denotes the volume of $S$. It is well known that for $1<p\le\infty$, we have
\[
\|Mf\|_p\le A_{K,p}\|f\|_p.
\]
Here we deal with the a related question: for what $K$ and $p$ does there exist some $\epsilon(K,p)>0$ such that
\begin{equation}\label{ineq:Lerner}
\|Mf\|_p\ge (1+\epsilon(K,p))\|f\|_p?
\end{equation}
This question was asked for the uncentered maximal operator with $K$ a ball by Lerner in \cite{Lerner2010}  and answered affirmatively for all $p<\infty$ in $\cite{IvanisviliJayeNazarov2017}$. In fact they showed this for general centrally symmetric convex $K$ and for uncentered maximal function defined by taking
\[
M_u f(x)=\sup_{\lambda>0, S=y+\lambda K\ni x}\frac{1}{|S|}\int_{S}f.
\]
 Similar positive results have been obtained for dyadic maximal functions \cite{Melas2017}, maximal functions defined over {\em $\lambda$-dense family} of sets, and {\em almost centered} maximal functions (see \cite{IvanisviliJayeNazarov2017} for details).
 
 This is closely related to the question of whether nonzero fixed points of $M$ exist in $L^p$; in fact if \eqref{ineq:Lerner} is satisfied, then no fixed points will exist. Korry \cite{Korry2001} proved that the centered maximal operator for $K$ a ball has fixed points if and only if $d \geq 3$ and  $p>\frac{d}{d-2}$, but a lack of fixed points does not imply that \eqref{ineq:Lerner} holds. On the other hand, by comparing $Mf(x) \geq C(d) M_{u}f(x)$, and using the fact that $\| M_{u}f\|_{L^{p}(\mathbb{R}^{n})} \geq (1+\frac{B(d)}{p-1})^{1/p}\|f\|_{L^{p}(\mathbb{R}^{n})}$ (see \cite{IvanisviliJayeNazarov2017}), one can easily conclude that (\ref{ineq:Lerner}) holds true whenever $p$ is sufficiently close to $1$. It is natural to ask what is the maximal $p_{0}(K)$ for which  if $1<p<p_{0}(K)$ then (\ref{ineq:Lerner}) holds. Similarly, we can ask what is the least $q_{0}(K)$ for which  if $q_0<p$ then $M$ has a fixed point in $L^p$.
 
 In  \cite{IvanisviliZbarsky2019}, the authors work in $d=1$ and prove $\eqref{ineq:Lerner}$ for all $p<1.5$, and for all $p$ for some special classes of functions $f$. Here, we prove it for all $p<\infty$ when $d=1$ or $d=2$, and also address the question of fixed points in dimension $d>2$, showing that for a generic shape $K$, $q_{0}(K)>\frac{d}{d-2}=q_0(B(0,1))$.
  
We first make a simplifying assumption. Since applying linear transformations to the shape does not change any of the inequalities we will be interested in, we will assume that
\begin{equation}\label{cond:superharmonic}
\int_K x_ix_j=\delta_{ij},
\end{equation}
which we can do by defining an  inner product $\<e_i,e_j\>=\int_K x_ix_j$, finding an orthonormal basis with respect to this inner product, and transforming $K$ by the change of basis matrix.

We prove the following:
\begin{theorem}\label{thm:growth}
For $d=1$ or $d=2$ and $1<p<\infty$, there exists $\epsilon=\epsilon(K,p)>0$ such that
\[
\|Mf\|_p\ge (1+\epsilon)\|f\|_p
\]
\end{theorem}
\begin{remark}
The constant $\epsilon(K,p)$ given by the proof can be computed for a given $K$ and $p$, but we cannot extract asymptotics in $p$.  It seems likely that the constant $\epsilon$ depends only on $p$ and not $K$, but it is not clear how to show this.
\end{remark}
As a lemma in the proof, we need the first part of the following theorem:
\begin{theorem}\label{thm:fixedpoints}
The following hold:
\begin{enumerate}
\item If $d=1$ or $d=2$, the only fixed points for $M$ in $L^\infty$ are constant functions.
\item
If $d\ge 3$, then there are no fixed points in $L^p$ for $p\le\frac{d}{d-2}$. Also, given $d$, there are coefficients $a_{ijk\ell}$ (given by fourth derivatives of the Green's function for Laplace's equation at a given point) such that if the shape $K$ satisfies
\[
\int_{K} \sum_{1\le i,j,k,\ell\le d} a_{ijk\ell}x^ix^jx^kx^\ell \ne 0
\]
then there is some $q>\frac{d}{d-2}$ (depending on $K$) such that there are no fixed points in $L^p$ for $p\le q$.
\end{enumerate}
\end{theorem}
The condition in the statement above is generic, so for generic shapes $K$, we have  $q_{0}(K)>q_0(B(0,1))$. We have also verified that the condition holds for the cube and the cross-polytope (the dual polytope of the cube). We also have a family of fourth-moment conditions obtained by rotating $K$, as well as analogous families of conditions for all even order moments (coming from the $2k$ order term in the Taylor expansion of the Green's function of Laplace's equation). We will not write out the conditions for higher order moments of $K$ explicitly, however one can extract them from the proof. If we think of the boundary of the polytope as given by the graph of some function $r(\theta)$ over $S^{n-1}$, the order $2k$ conditions reduce to some (possibly all) of the order $2k$ spherical harmonic coefficients of $r^{2k+1}$ being 0. It is not clear how easy it is to satisfy all these conditions simultaneously.
\begin{question}
If the maximal function for a convex centrally symmetric shape $K$ has fixed points in $L^p$ for all $p>\frac{d}{d-2}$, is $K$ necessarily a ball?
\end{question}
\section{Dimensions 1 and 2}\label{sec:d12}
We look at two closely related notions in this section. First, we take the increasing sequence of functions $\{ M^n\mathbbm{1}_{B(0,1)}\}$ and define the pointwise limit $f=\lim_{n\to\infty} M^n\mathbbm{1}_{B(0,1)}$.  Second, we look at fixed points of $M$. These are related because it is easy to see that $f$ is a fixed point of $M$ (for this we need that $\|f\|_\infty=1$) and conversely, if $g$ is a nonzero fixed point of $M$, then $g$ must be greater than $\delta>0$ on some disk, so up to rescaling,
\[
g=\lim_{n\to\infty} M^ng\ge \lim_{n\to\infty} M^n\mathbbm{1}_{B(0,1)}=f.
\]
Thus, understanding the function $f$ is intimately tied to understanding what function spaces contain fixed points of $M$.

The following argument for dimensions 1 and 2 is the argument from $\cite{Korry2001}$, generalized to shapes $K$ other than the disk when $d=2$. For shapes other than the disk, we also use the idea of mollification and using Taylor series from $\cite{IvanisviliJayeNazarov2017}$. Since the authors of that paper are looking at off-center maximal functions, they only need the linear term of the Taylor expansion; we use the quadratic term in this section and the quartic and higher-order terms in Section~\ref{sec:dbig}.

Now suppose that we have an $L^\infty$ fixed point $g$ of $M$ in $d=1$ or $d=2$.  We want to prove that $g$ is constant. We let $\tilde g=g*\eta$ for $\eta$ a standard mollifier. Then
\[
M\tilde g\le (Mg)*\eta=g*\eta=\tilde g
\]
so $\tilde g$ is also a fixed point of $M$ and is smooth. Apply the definition of $M$ (see \eqref{Mdef}) for $\lambda$ small to get
\[
M\tilde g(x)\ge\tilde g(x)+\frac{\lambda^2}{2}\Delta \tilde g(x)+O(\lambda^3).
\] 
We get no $\lambda$ term because we are expanding around the center of mass, and the $\lambda^2$ term is as it is because of condition~$\eqref{cond:superharmonic}$. If $\Delta \tilde g(x)>0$, we then have that $M\tilde g(x)>\tilde g(x)$, which contradicts $\tilde g$ being a fixed point. Thus $\tilde g$ is superharmonic.

$\tilde g$ being constant now follows from the arguments in the classification of fixed points in $\cite{Korry2001}$ since superharmonicity was the only fact about fixed points that was used in that paper (In that paper, Korry's argument forces $L^\infty$ superharmonic functions to be constant in $d=1$ and $d=2$). Since any mollification of $g$ is constant, $g$ in turn is constant.

We will now prove Theorem~\ref{thm:growth}. We follow a similar argument to $\cite{IvanisviliJayeNazarov2017}$. The constants are chosen with the $d=2$ case in mind; however the argument as written also covers the $d=1$ case. Since $M$ is Lipschitz on $L^p$ and since functions are approximated arbitrarily well in $L^p$ by continuous compactly supported functions, we assume without loss of generality that $f$ is continuous and compactly supported. Take some small $\delta_1>0$ to be chosen later.

Given some $\mu>0$, take all shapes $\{S_i=x_i+\lambda_i K\}$ on which the average of $f$ is exactly equal to than $\mu(1-\delta_1)$. Note that if $f(x)\ge\mu$, then by continuity, there is some $\epsilon>0$ such that the average of $f$ over $x+\epsilon K$ is greater than  $\mu(1-\delta_1)$. Also, as $R\to\infty$, the average of $f$ over $x+R K$ goes to 0. Thus by the intermediate value theorem, $x$ is the center of some $S_i$.

Note that since $K$ is convex, $K+\delta_1 K\subseteq (1+\delta_1)K$. By using this fact and taking slight rescalings and shifts of $S_i$, we get that $Mf\ge \mu(1-\delta_1)/(1+\delta_1)^2$ on $x_i+\delta_1 \lambda_i K$. We then have that
\[
\lim_{n\to\infty} M^n \mathbbm{1}_{\delta_1 K}
\]
is a constant function (by Theorem 1), so
\[
\lim_{n\to\infty} M^n \mathbbm{1}_{\delta_1 K}=1,
\]
so by choosing $n$ large enough, we have that
\[
 M^{n-1} \mathbbm{1}_{\delta_1 K}\ge 1/(1+\delta_1)
\]
on all of $2K$ excluding a set of volume $\delta_1/|K|$. Applying $M$ again at every point of $K$ and using $\lambda=1$, we get that
\[
 M^n \mathbbm{1}_{\delta_1 K}\ge (1-\delta_1)/(1+\delta_1)
\]
on $K$. Note that this $n$ is independent of the function $f$, but depends on the shape $K$. Then
\begin{equation}\label{manyiter}
M^{n+1}f\ge M^n\left(\frac{\mu(1-\delta_1)}{(1+\delta_1)^2}\mathbbm{1}_{x_i+\delta_1\lambda_i K}\right)\ge  \mu(1-\delta_1)^2/(1+\delta_1)^3 \mathbbm{1}_{S_i}.
\end{equation}

Now we use the  Besicovitch covering lemma as given in Appendix~\ref{appendix} to extract a countable subset of $\{S_{i_j}\}_{j\in\N}$. We take $E=f(x)\ge\mu$  and recall that each $x\in E$ was the center of some $S_i$, in particular $x\in \bigcup S_{i_j}$. Also, each point is covered at most $B(d)$ times where $B(d)$ is the constant in the statement of the Besicovitch covering lemma. Thus $\sum |S_{i_j}|$ is finite, which we will use for rearranging sums. We let $\alpha(x)=|\{j\mid x\in S_{i_j}\}|$. Note that on $ \bigcup S_{i_j}$, we have that $1\le\alpha\le B(d)$. We have that
\[
0=\sum_j\int_{S_{i_j}}(f-\mu(1-\delta_1))=\int_{\bigcup S_{i_j}}\alpha(x)(f(x)-\mu(1-\delta_1))dx\ge \mu|\{f\ge 2\mu\}|-B(d)\mu\left|\{f<\mu\}\cap\left(\bigcup S_{i_j}\right)\right|.
\]
Thus
\[
|\{M^{n+1}f\ge  \mu(1-\delta_1)^2/(1+\delta_1)^3\}|\ge \left|\bigcup S_{i_j}\right|\ge |\{f\ge \mu\}|+|\{f\ge 2\mu\}|/B(d)
\]
where we used \eqref{manyiter} for the first inequality.

Multiplying both sides of the above inequality by $p\mu^{p-1}$ and integrating $\mu=0$ to $\infty$, we get
\[
(1+\delta_1)^{3p}/(1-\delta_1)^{2p}\|M^{n+1}f\|_p^p\ge \|f\|_p^p+\|f\|_p^p/(2^pB(d))
\]
so by picking $\delta_1$ sufficiently small, we get that
\begin{equation}\label{ineq:stronggrowth}
\|M^{n+1}f\|_p^p\ge (1+\epsilon(K,p))\|f\|_p^p.
\end{equation}

We now copy the argument from \cite{IvanisviliZbarsky2019}. Suppose that $\|Mf-f\|_p<\tilde\epsilon\|f\|_p$ for some $\tilde\epsilon$ to be chosen later. From the subadditivity of the maximal operator, it follows that $\|M\phi_1-M\phi_2\|_p\le A_{K,p}\|\phi_1-\phi_2\|_p$, so
\[
\|M^{n+1}f-f\|_p\le \sum_{j=1}^{n+1}\|M^jf- M^{j-1}f\|_p\le\sum_{j=1}^{n+1}A_{K,p}^{j-1}\|Mf- f\|_p<\left(\tilde\epsilon\sum_{j=1}^{n+1}A_{K,p}^{j-1}\right)\|f\|_p
\]
which contradicts $\eqref{ineq:stronggrowth}$ for $\tilde\epsilon=\tilde\epsilon(p)$ sufficiently small. Thus $\|Mf-f\|_p\ge\tilde\epsilon\|f\|_p$, so
\[
\|Mf\|_p^p=\int (Mf)^p\ge\int f^p+(Mf-f)^p=\|f\|_p^p+\|Mf-f\|_p^p\ge\left(1+\tilde\epsilon^p\right)\|f\|_p^p,
\]
which proves Theorem~\ref{thm:growth}.

The constant $\epsilon(K,p)$ given by the proof can be computed for a given $K$ and $p$. However, this involves getting a bound on $n$, which involves understanding how fast $\{M^n \mathbbm{1}_{\delta_1 K}\}$ converges to 1. For a given $K$ and $p$, this can be done by approximately calculating the sequence $\{M^n \mathbbm{1}_{\delta_1 K}\}$, but this is inefficient and does not give asymptotics in $p$. It seems likely that the constant $\epsilon$ depends only on $p$ and not $K$, but it is not clear how to show this.
\section{$d\ge 3$}\label{sec:dbig}
In this section, we prove the $d\ge 3$ case of Theorem~\ref{thm:fixedpoints}. In order to prove that there is no fixed point of $M$ in $L^p$, we will show that
\begin{equation}\label{ineq:notinLp}
\lim_{n\to\infty} M^n\mathbbm{1}_{B(0,1)}\ge 2^{-d/p}\mathbbm{1}_{B(0,2)}.
\end{equation}
Iterating \eqref{ineq:notinLp}, we get that $\lim_{n\to\infty} M^n\mathbbm{1}_{B(0,1)}$ is not in $L^p$. Since any fixed point of $M$ lies above some rescaling of $\mathbbm{1}_B(0,1)$, we get that the fixed point is not in $L^p$.

We now turn to understanding for what $p$ we can prove \eqref{ineq:notinLp}. Suppose that we have a fixed point $g\in L^p(\R^d)$ for $d\ge 3$.  We let $\tilde g=g*\eta$ for $\eta$ a standard mollifier. Then
\[
M\tilde g\le (Mg)*\eta=g*\eta=\tilde g
\]
so $\tilde g$ is also a fixed point of $M$ and $\tilde g\in L^p$. Apply the definition of M for $\lambda$ small to get
\[
M\tilde g(x)\ge\tilde g(x)+\frac{\lambda^2}{2}\Delta \tilde g(x)+O(\lambda^3).
\] 
We get no $\lambda$ term because we are expanding around the center of mass, and the $\lambda^2$ term is as it is because of condition~$\eqref{cond:superharmonic}$. If $\Delta \tilde g(x)>0$, we then have that $M\tilde g(x)>\tilde g(x)$, which contradicts $\tilde g$ being a fixed point. Thus $\tilde g$ is superharmonic.

Since $g$ yields a superharmonic function after any mollification, we can see from the weak formulation of superharmonic functions that it in turn is superharmonic. Now take
\[
g=\lim_{n\to\infty} M^n\mathbbm{1}_{B(0,1)}.
\]
Then $g$ is a fixed point of $M$, thus superharmonic.
We take some large radius $R$. When $|x|=R$, we have $g\ge 0$. When $|x|=1$, we have $g\ge 1$. Solving Laplace's equation, we get that on the annulus $1<|x|<R$, we have $g\ge \frac{|x|^{2-d}-R^{2-d}}{1-R^{2-d}}$.
Taking $R\to \infty$, we get $g\ge |x|^{2-d}$. This already proves  $\eqref{ineq:notinLp}$ for $p\le \frac{d}{d-2}$. We now let $h(x)=|x|^{2-d}$ and investigate $Mh(x)$ on $\{|x|=3\}$.

By applying the definition of M for $\lambda$ small and the Taylor expansion of $h$, we get
\begin{equation}\label{expr:moment}
Mh(x)=h(x)+\frac{\lambda^4 }{24|K|}\sum_i\sum_j\sum_k\sum_\ell\partial_{ijk\ell}h\int_K  x_ix_jx_kx_\ell+O(\lambda^5).
\end{equation}
where the quadratic term vanishes because $\Delta h=0$ and $K$ satisfies $\eqref{cond:superharmonic}$ and the linear and cubic terms vanish because $K$ is centrally symmetric.
When the quartic term of $\eqref{expr:moment}$ is positive, we can take $\lambda$ small, and get that $Mh(3,0,\ldots,0)>3^{2-d}$. By continuity, this is true on a small region around this point. Since $g\ge h$ is a fixed point of $M$ and superharmonic, a lower bound on $g$ is the solution to Laplace's equation on the annulus $1<|x|<3$ with boundary data given by $1$ on $\{|x|=1\}$ and $Mh$ on $\{|x|=3\}$. Since $Mh>h$ on some portion of the boundary, we have that this harmonic function is strictly greater than $h$ everywhere in the interior. By compactness of $\overline{B(0,2)}$, we then get
\[
\lim_{n\to\infty} M^n\mathbbm{1}_{B(0,1)}\ge 2^{-d/p}\mathbbm{1}_{B(0,2)}
\]
for some $p>d/(d-2)$.

When the quartic term of $\eqref{expr:moment}$ is negative, we note that
\[
\int_{S^{n-1}}\dashint_{3z+\epsilon K} h  dz=3^{2-d}
\]
because averaging over points on the sphere is the same as averaging over all possible rotations of $K$ at one point, and the latter will give $h(3,\ldots,0)$ by the mean value property of harmonic functions. Thus if the quartic term of $\eqref{expr:moment}$ is negative at some point, there will be some other points on the sphere $\partial B(0,3)$ such that $Mh>3^{2-d}$ in a neighborhood of that point. Then by the same argument as above, we get
\[
\lim_{n\to\infty} M^n\mathbbm{1}_{B(0,1)}\ge 2^{-d/p}\mathbbm{1}_{B(0,2)}
\]
for some $p>d/(d-2)$. Thus we can only hope to avoid having such a $p>d/(d-2)$ if the quartic term of $\eqref{expr:moment}$ is 0. This is precisely the condition in the statement of Theorem~\ref{thm:fixedpoints}, completing the proof of that theorem.

To avoid having some $p>d/(d-2)$ work, the quartic term of \eqref{expr:moment} must be 0. Also, all the expressions we get from it by rotating the coordinate system must be 0. When these expressions are 0, we can expand the Taylor series of $h$ to 6th order and get conditions that some combinations of the 6th moments must also be 0, and so on. This is the infinite family of conditions alluded to below the statement of Theorem~\ref{thm:fixedpoints}.

\section{Acknowledgements}
The author would like to thank Paata Ivanisvili for earlier discussions related to these problems. This material is based upon work supported by the National Science Foundation Graduate Research Fellowship Program under Grant No. DGE-1656466. Any opinions, findings, and conclusions or recommendations expressed in this material are those of the author and do not necessarily reflect the views of the National Science Foundation.

\appendix
\section{Besicovitch Covering Lemma}\label{appendix}
\begin{lemma}
Let $K\subset\R^d$ be convex, compact, and centrally symmetric. Suppose we have a bounded set $E$ and some constant $\Lambda>0$, and take
\[
A=\{x+\lambda_x K\}_{x\in E}
\]
where $0<\lambda_x<\Lambda$ for each $x\in E$. Then there is some countable $\tilde A\subset A$ so that the sets $\tilde A$ cover all of $E$, and each point is covered at most $B(d)$ times, where $B(d)$ is a constant depending only on $d$.
\end{lemma}
Note that this formulation is the same as the formulation for balls in some norm on $\R^d$. A proof can be obtained by a straightforward modification of the standard proof of the Bescovitch covering lemma for the usual norm on $\R^n$, which can be found for instance in \cite{PerttiGeometryOfSetsAndMeasures}. The appendix of \cite{IvanisviliJayeNazarov2017} also gives references which they claim can be found to contain the proof of the Besicovitch Covering Lemma for arbitrary norm with enough digging.
\bibliographystyle{abbrv}
\bibliography{zbarskybib}
\end{document}